\documentclass[11pt]{article}

\usepackage{amsmath,amsthm,amscd,amssymb}
\usepackage[alphabetic]{amsrefs}
\usepackage[margin=1.25in]{geometry}
\usepackage[colorlinks, citecolor = blue]{hyperref}
\usepackage[capitalize]{cleveref}
\usepackage{multirow}
\usepackage{verbatim}

\newcommand{\Q}{\mathbb Q}
\newcommand{\Z}{\mathbb Z}

\newcommand{\calO}{\mathcal O}

\newcommand{\frakp}{\mathfrak p}

\newcommand{\frako}{\mathfrak o}

\newcommand{\bmx}{\begin{pmatrix}}
\newcommand{\emx}{\end{pmatrix}}
\newcommand{\disc}{\text{disc}}
\newcommand{\lev}{\text{level}}
\newcommand{\la}{\langle}
\newcommand{\ra}{\rangle}
\newcommand{\ds}{\displaystyle}

\newcommand{\quat}[3]{\left(\frac{#1,#2}{#3}\right)}
\newcommand{\hs}[2]{\left(\frac{#1}{#2}\right)}
\newcommand{\qam}[1]{\left\{\bmx \alpha & b\beta\\ \bar{\beta} & \bar{\alpha}\\ \emx: \alpha,\beta\in #1\right\}}
\newcommand{\qom}[3]{\left\{\bmx \alpha & #1\beta\\ \bar{\beta} & \bar{\alpha}\\ \emx: \alpha\in #2,\beta\in #3\right\}}

\renewcommand{\mod}{\, \, \mathrm{mod} \, \,}

\newtheorem{lem}{Lemma}
\newtheorem{prop}[lem]{Proposition}
\newtheorem{thm}[lem]{Theorem}

\newtheorem{ex}[lem]{Example}

\numberwithin{lem}{section}
\numberwithin{equation}{section}

\begin{document}

\title{Constructing non-maximal orders in quaternion algebras}
\author{Jordan Wiebe}
\date{\today}
\maketitle

\begin{abstract}
We present an explicit basis for orders of arbitrary level $N>1$ in definite rational quaternion algebras.  These orders have applications to computations of spaces of elliptic and quaternionic modular forms.
\end{abstract}

\section{Introduction}
Let $S_k(N)$ denote the space of elliptic cusp forms of weight $k$ on $\Gamma_0(N)$ with trivial character.  Denote by $B$ a definite quaternion algebra over $\Q$, and denote by $\calO$ an order in $B$.  In 1940, Hecke conjectured that for a prime $p$, a basis for $S_2(p)$ could be obtained via the theta series associated to a set of $\calO$-ideal class representatives, where $\calO$ is a maximal order in the definite quaternion algebra ramified at $p$ and $\infty$.  In 1956, Eichler \cite{eichler-56} proved that there was indeed a basis for $S_2(p)$ taken from a more general collection of theta series associated to $\calO$ obtained via certain arithmetically defined matrices associated to the order, called Brandt matrices.  Hijikata, Pizer and Shemanske \cite{hps-memoirs} generalized Eichler's work to arbitrary level in 1989 by providing orders of level $N=p^rM$ in a definite quaternion algebra ramified at $p$ and $\infty$.  More recently, Martin \cite{martin-basis} treated the basis problem using orders in algebras with more general discriminant.

Much of the literature on the subject of these orders involves fixing a definite quaternion algebra $B$ ramified at $p$ and $\infty$, and explicit bases for the orders presented are quite limited.  In particular, Pizer \cite{pizer-algorithm} presents bases for maximal orders of the definite algebras ramified at a single finite prime.  Albert \cite{albert} and Ibukiyama \cite{ibuk} also present bases for maximal orders of definite quaternion algebras.  Pacetti and Rodriguez-Villegas \cite{p2-orders} also provide bases for orders of level $p^2$.  More generally, orders of level $N$ can be used to construct modular forms of level $N$.  Note that there are other approaches to this, for instance modular symbols or an approach of Demb\'el\'e \cite{dembele} which only requires the use of maximal orders.  However, explicit construction of non-maximal orders is also useful for computations of quaternionic modular forms via Brandt matrices, which we plan to address in the future.  Specifically, computations of these quaternionic modular forms will allow us to address questions raised in \cite{martin-congruence}.

We present here an explicit basis for orders of arbitrary admissible level $N>1$ in definite rational quaternion algebras, where an admissible level $N$ is one in which the discriminant of the quaternion algebra divides $N$ (note that this is a necessary condition to obtain an order of level $N$ in $B$).  These results have been checked for $\Delta\leq 1000$ and $N\leq 10,000$ in Sage.  Furthermore, we present these orders in arbitrary definite rational quaternion algebras in every case except the case where $v_2(N)=2$ and the discriminant of $B$ is even.  In this case, we can construct an order with level $N=4N'$ (with $N'$ odd) in a quaternion algebra with even discriminant if $\prod_{p\mid RM_1, p\neq 2} p\equiv 1\mod 4$, where we have written our level as $N=RM$, where the discriminant of $B$ divides $R$, and split $M$ into $M_1$ and $M_2$, where the factors of $M_1$ are the primes in $N$ which have odd exponent, and the factors of $M_2$ are the primes in $N$ which have even exponent.

Our result uses a careful choice of the presentation of our definite quaternion algebra $B=\quat{a}{b}{\Q}$ via $a,b\in\Z$, allowing for computation of a space of modular forms of weight $2k$ and level $N$.  When we combine our result with twists of forms of lower level, we can obtain the entire newspace of weight $2k$ and level $N$ (see \cite{hps-orders}, \cite{martin-basis}).  Moreover, our result can be used to compute quaternionic modular forms via Brandt matrices.  Note that our result is for definite rational quaternion algebras, but a nearly identical argument (with different conditions $\mod 8$) will work for indefinite quaternion algebras.  The general description of our basis is somewhat complex, so for simplicity we state an explicit basis for the special case with odd level $p^r$:

\begin{thm}\label{thm-1-1}

Let $B$ be a definite quaternion algebra with discriminant $\Delta_B=p$ odd.  Take \[
	a,b = \begin{cases}
		-q, -p & \text{if $r$ is odd and $p\equiv 1\mod 4$, with $q$ nonsquare $\mod p$ and $3\mod 4$}\\
		-p, -q & \text{if $r$ is odd and $p\equiv 3\mod 4$, with $q$ square $\mod p$ and $1\mod 8$}\\
		-qp, -p & \text{if $r$ is even and $p\equiv 1\mod 4$, with $q$ nonsquare $\mod p$ and $3\mod 4$}\\
		-p, -q & \text{if $r$ is even and $p\equiv 3\mod 4$, with $q$ square $\mod p$ and $1\mod 8$}\\
	\end{cases}
\]
Then we can represent $B$ as $\quat{a}{b}{\Q}$.  Furthermore, put $f = p^n$ if $r$ is odd, and $f = p^{n-1-v_p(b)}$ if $r$ is even, and select $x$ with $x^2\equiv -p \mod q$, and let $z$ be given by the Euclidean algorithm for finding $y(-q)+z(2x)=1$, $u$ be given by using the Euclidean algorithm to write $v(q)+w(2x)=1$, and setting $0\leq u < 2q$ such that $u\equiv vq + 2w \mod 2q$, and $z'$ is given by choosing $0\leq z' < 2q$ with $z'\equiv 4z \mod 2q$.  Then the order \[
	\calO = \begin{cases}
		\Z\left\la \frac{q+i+2zk}{2q}, \frac{2i+z'k}{2q}, \frac{f(j+k)}{2}, fk \right\ra & \text{if } p\equiv 1\mod 4\\
		\Z\left\la 1, \frac{1+i}{2}, \frac{f(j+uk)}{2q}, fk \right\ra & \text{if } p\equiv 3\mod 4\\
	\end{cases}
\]
has level $N=p^{r}$.

\end{thm}

The above theorem splits into two cases, one where $r$ is odd and the other where $r$ is even.  In the second case, $p$ is ramified in $K$, making the construction more complicated, as well as the structure of the space of associated theta series (see \cite{hps-memoirs}, \cite{martin-basis}).  Our result is stated in full generality in Theorem \ref{big-thm}.

In Section 2, we will cover preliminaries on quaternion algebras.  In Section 3, we will embed our quaternion algebra $B$ in $M_2(K)$, for $K$ a quadratic field, and examine its level locally in cases based on the splitting/ramification of $K_p$ and $B_p$.  In Section 4, we will construct a global order $\calO$ of $B$ using the local results from the previous section, with level $N\cdot q$ for a suitable auxiliary prime $q$.  We will construct this order in cases, based on the behavior of $2$ in the quadratic field $K$ and in the algebra $B$.  We will also calculate the basis for this order.  In Section 5, using a technique of Voight \cite{voight-max} we will lower the level of our order constructed in the previous section from $N\cdot q$ to $N$, and calculate the new basis for this order.  In Section 6, we will present our general result, as well as a few special cases, including an order of level $p^r$ for the algebra ramified at a single prime $p$.  Finally, in Section 7 we will present examples of using our construction of orders to compute spaces of modular forms.

\textbf{Acknowlegements:} I would like to thank Kimball Martin for his guidance and insight, generously shared with me throughout the process of this work.  I am also grateful to John Voight for helpful discussions.

\section{Preliminaries on Quaternion Algebras}
Recall that a quaternion algebra over $\Q$ is a four-dimensional central simple $\Q$-algebra.  Note that any quaternion algebra over $\Q$ is either a noncommutative division algebra or the split matrix algebra $M_2(\Q)$.  We can construct quaternion algebras using the Hilbert symbol $\quat{a}{b}{\Q}$ to denote the quaternion algebra with $\Q$-basis $1,i,j,k$ and multiplication satisfying \[
	i^2=a, j^2=b, \text{ and } ij=-ji=k.
\]
Indeed, any quaternion algebra can be constructed in this way for some $a,b\in\Z$.

The splitting behavior of our quaternion algebra is described as follows:

\begin{lem}\label{lem-quatsplit}
Suppose that $p$ is an odd rational prime and $a,b\in\Z$ are nonzero and square-free.  Assume $v_p(a)\leq v_p(b)$.  Then $\left(\frac{a,b}{\Q_p}\right)$ is division (i.e., ramified) if and only if 
\begin{enumerate}
	\item $p\nmid a$, $p\mid b$, and $a$ is a nonsquare mod $p$; or 
	\item  $p\mid a$, $p\mid b$, and $-a^{-1}b$ is a nonsquare mod $p$.
\end{enumerate}

Alternatively, if $p=2$, then we have \[
	(a,2)_{\Q_2}=\hs{a}{2}=\begin{cases}
	+1 & \text{if } a\equiv 1, 7 \text{ mod } 8\\
	-1 & \text{if } a\equiv 3, 5 \text{ mod } 8
\end{cases}.
\]
Furthermore, if $a$ and $b$ are odd primes we have $(a,b)_{\Q_2}=(-1)^{\frac{a-1}{2}\frac{b-1}{2}}$.
\end{lem}

This lemma follows from known calculations of Hilbert symbols over $\Q_p$.

There are $3$ possibilities for the behavior of $K_p$, and 2 for $B_p$, so we have the means of describing the behavior of our quaternion algebra $B$ with $K=\Q(\sqrt{a})\subset B$ in six cases.  However, in order for our quadratic field $K$ to be contained in $B$, $B$ must not ramify if $K$ is split; in other words, we can omit one of the cases: \[
	\renewcommand{\arraystretch}{1.3}
	\begin{tabular}{| c | c | c |}
		\hline
		$K_p$ split, $B_p$ split & $K_p$ ramified, $B_p$ split & $K_p$ unramified, $B_p$ split\\
		\hline
		$\times$ & $K_p$ ramified, $B_p$ ramified & $K_p$ unramified, $B_p$ ramified\\
		\hline
	\end{tabular}
\]
	
We also note here that the Hilbert symbols defined above have many helpful properties, including $(a,b)_F\cdot (a,c)_F = (a,bc)_F$ if $F$ is $p$-adic.  This will prove useful for calculating the behavior of $B_2$, the localization of our algebra $B$ at $2$.

\section{Local Orders}\label{local-orders}
Recall that we can embed our quaternion algebra $B=\quat{a}{b}{\Q}$ in $M_2(K)$, where $K=\Q(\sqrt{a})$ as \[
	B=\qam{K=\Q(\sqrt{a})}.
\]
The form of an order varies depending on whether $K$ is split, ramified, or unramified, and also varies based on whether $B$ is split or ramified.  Therefore, we will examine orders in each case separately.  In particular, we may consider the local algebra \[
	B_p = B\otimes_{\Q} \Q_p = \qam{K_p=K\otimes\Q_p}
\]
at each prime, and examine orders locally, split into the cases above for the behavior of $K_p$ and $B_p$.  Consider the following orders, referred to as residually inert orders by Voight: for finite primes that ramify in $B$, a special order $\calO_p$ of $B_p$ has a quadratic extension $K_p$ of $\Q_p$ and a positive integer $v(p)$ (odd if $K_p$ is unramified) so that $\calO_p=\frako_{K_p}+\mathfrak{P}^{v-1}_{B_p}$, where $\mathfrak{P}_{B_p}$ is the unique maximal ideal of the unique maximal order $\calO_{B_p}$ of $B_p$.  Note that these orders were called special orders by Hijikata, Pizer, and Shemanske in \cite{hps-orders}.  For our purposes, we will construct $\calO$ to be a residually inert order for primes which ramify in $B$ and in $K$ (using $K$ as our quadratic field), and to be Eichler (residually split) for primes which split in $B$.

For the remainder of this section, fix $a,b\in\Z$ square-free and coprime.

\subsection{$B_p$ is split}
Assume that $B_p$ is split.  Then the standard Eichler order $\calO_{B_p}(n)$ of level $n=p^k$ has the form \[
	\calO_B(n) = \bmx \Z_p & \Z_p \\ p^k\Z_p & \Z_p \\ \emx,
\]
and all Eichler orders of level $n$ are conjugate to $\calO_B(n)$.

\subsubsection{$K_p$ is split}\label{split-split}
If $K_p$ is split, we have $K_p=\Q_p\oplus\Q_p$ and $\frako_{K_p}=\Z_p\oplus\Z_p$.  Consider the order \[
	\calO_p = \qom{b}{\frako_{K_p}}{f\frako_{K_p}} = 
	\calO_p = \left\{\bmx (x,y) & bf(z,w) \\ f(w,z) & (y,x) \\ \emx: x,y,z,w\in\Z_p\right\}.
\]
We now conjugate and simplify: \[
\bmx
	f & \\
	 & 1\\
\emx
\bmx
	(x,y) & bf(z,w)\\
	f(w,z) & (y,x)\\
\emx
\bmx
	\frac{1}{f} & \\
	 & 1\\
\emx
=
\bmx
	(x,y) & bf^2(z,w)\\
	(w,z) & (y,x)\\
\emx.
\]
So we can identify these matrices with pairs of matrices $\left(\bmx x & bf^2z\\ 		w & y\\ \emx,\bmx y & bf^2w\\ z & x\\ \emx\right)$, and we have \[
	\calO_p \simeq \left\{\bmx
	x & bf^2z\\
	w & y\\
	\emx: x,y,z,w\in\Z_p\right\}.
\]
This is an Eichler order of level $p^{2v_p(f)+v_p(b)}$ in $B_p=M_2(\Q_p)$.

\subsubsection{$K_p$ is ramified}\label{split-ram}
Now assume that $K_p$ is ramified.  Note that $\frako_{K_p}=\Z_p[\sqrt{a}]$.
	
If $B_p$ is split, then $b$ must be a norm from $\frako_{K_p}^\times$, so we can write $b=u\bar{u}$ for some $u\in\frako_{K_p}^\times$.  Making the substitution $\beta\mapsto \bar{u}^{-1}\beta$ gives us \[
	B_p = \left\{\bmx
	\alpha & u\beta\\
	u^{-1}\bar{\beta} & \bar{\alpha}\\
\emx: \alpha,\beta\in K_p\right\}.
\]
Consider the order \[
	\calO_p = \left\{\bmx
	\alpha & u\beta\\
	u^{-1}\bar{\beta} & \bar{\alpha}\\
\emx: \alpha\in\Z_p+g\frako_{K_p},\beta\in f(\Z_p+g\frako_{K_p})\right\}.
\]
Now let $\ell=\bmx \sqrt{a} & -\sqrt{a}\\ 1 & 1\\ \emx \bmx u^{-1} & \\  & 1\\ \emx			$, and write $\alpha = x+gy+gz\sqrt{a}$ and $\beta=fp+fgq+fgr\sqrt{a}$.  Then the conjugation $\ell\calO_p\ell^{-1}$ gives us \[
	\ell\bmx \alpha & u\beta \\ u^{-1}\bar{\beta} & \bar{\alpha}\\ \emx \ell^{-1}=
	\frac{1}{2u^{-1}\sqrt{a}}\cdot\bmx u^{-1}\sqrt{a}[(\alpha+\bar{\alpha})-(\beta+\bar{\beta})] & u^{-1}a[(\alpha-\bar{\alpha})+(\beta-\bar{\beta})] \\ u^{-1}[(\alpha-\bar{\alpha})-(\beta-\bar{\beta})] & u^{-1}\sqrt{a}[(\alpha+\bar{\alpha})+(\beta+\bar{\beta})] \\ \emx.
\]
Now if $\alpha\in\Z_p+g\frako_{K_p}$ and $\beta\in f(\Z_p+g\frako_{K_p})$ then $\alpha+\bar{\alpha}=2x+2gy$, $\alpha-\bar{\alpha}=2gz\sqrt{a}$, $\beta+\bar{\beta}=2fp+2fgq$, and $\beta-\bar{\beta}=2fgr\sqrt{a}$.  This gives us \[
	\ell\calO_p\ell^{-1} = \bmx x+gy-fp-fgq & a(gz+fgr) \\ gz-fgr & x+gy+fp+fgq \\ \emx
\]
\[
	= \bmx \Z_p & p^{v_p(a)+v_p(g)}\Z_p \\ p^{v_p(g)}\Z_p & \Z_p \emx = \bmx \Z_p & p^{v_p(a)+2v_p(g)}\Z_p \\ \Z_p & \Z_p \emx.
\]
So we have an explicit conjugation of $B_p$ to $M_2(\Q_p)$ that clearly expresses $\calO_p$ as an Eichler order of level $p^{v_p(a)+2v_p(g)}$.

If $p=2$, we must be more careful, since it is possible for $2$ to be ramified in $K$ but for $2\nmid a$.  In particular, if $a\equiv 3 \mod 4$ but $2\nmid a$, then the order described above has $v_2(N)=2v_2(g)+2$.  Alternatively, if $p=2$ and $2\mid a$, then the order described above has $v_2(N)=2v_2(g)+8$.  Finally, if $p=2$ and $2\mid b$, then the order described above has $v_2(N)=2v_2(g)+1$.

\subsubsection{$K_p$ is unramified}
If $K_p$ is unramified, then $\frako_{K_p}=\Z_p[\sqrt{a}]$ unless $p=2$ and $a\equiv 1 \mod 4$, when we have $\frako_{K_2}=\Z_2[\frac{1+\sqrt{a}}{2}]$.  An identical argument as the previous section gives us \[
	\calO_p = \left\{\bmx
	\alpha & u\beta\\
	u^{-1}\bar{\beta} & \bar{\alpha}\\
\emx: \alpha\in\Z_p+g\frako_{K_p},\beta\in f(\Z_p+g\frako_{K_p})\right\}.
\]
with level $p^{v_p(a)+2v_p(g)}$.  In this case, since $K_p$ is unramified, we have $v_p(a)=0$, so our level is $p^{2v_p(g)}$.

If $p=2$ and $a\equiv 1 \mod 4$, then we have $\frako_{K_2}=\Z_2\left[\frac{1+\sqrt{a}}{2}\right]$.  The basis for our order $\calO_2$ is $\Z_2\la 1, \frac{1+i}{2}, fj, f\cdot\frac{j+k}{2}\ra$, and a quick calculation shows that $v_2(\disc(\calO))=2v_2(f)$.  So our level is $2^{2v_2(f)}$.

\subsection{$B_p$ is ramified}
Assume that $B_p$ is ramified.  We again break into cases based on whether $K_p$ is ramified or unramified.

\subsubsection{$K_p$ is ramified}\label{ram-ram}
Assume $K_p$ is ramified.  Here we use the residually inert orders defined previously.  We know that $\mathfrak{P}_{B_p}=\varpi_{B_p}\calO_{B_p}$, and $\mathfrak{P}_{B_p}^{v-1}=\{x\in\calO_{B_p}: N(x)\in\frakp^{v-1}\}$.  Now $\calO_{B_p}$ is the maximal order of $B_p$, a local division algebra, which we know is of the form \[
	\calO_{B_p}=\qom{b}{\frako_{K_p}}{\frako_{K_p}}
\]
if $p\neq 2$.  When $p=2$, \[
	\calO_{B_2} = \qom{2}{\frako_{K_2}}{\frako_{K_2}}
\]
is maximal when $K_2$ is unramified.

  We require that $p\mid b$ if $K_p$ is unramified, and $p\nmid b$ if $K_p$ is ramified to obtain the maximal orders.

Now consider an element $x\in \calO_B$: \[
	x = \bmx \alpha & b\beta\\ \bar{\beta} & \bar{\alpha} \emx = \bmx \alpha & \\  & \bar{\alpha} \emx +  \bmx  & b\beta\\ \bar{\beta} & \emx.
\]
For $x$ to be an element of our residually inert order $\calO_p$, we need $\bmx \alpha & \\ & \bar{\alpha} \emx\in \frako_{K_p}$, so $\alpha\in\frako_{K_p}$.  We also need $y=\bmx & b\beta\\ \bar{\beta} & \emx\in\mathfrak{P}_B^{v-1}$ to obtain level $p^{v}$.  Now $y\in\mathfrak{P}_B^{v-1}$ if and only if $N(y)\in p^{v-1}\Z_p$, and we know that $N(y)=-b\beta\bar{\beta}$.  So we know that we need $\beta\bar{\beta}\in p^{v-1-v_p(b)}\Z_p$.  Now we write $\beta=u\varpi_{K_p}^m$ for a uniformizer $\varpi$ and a unit $u$.  Now since $K_p$ is ramified, we choose $\varpi_{K_p} = \sqrt{vp}$.  So $\beta\bar{\beta} = u\bar{u}\cdot v^mp^m$ when $K_p$ is ramified.  So when $K_p$ is ramified, then if $\beta\in p^{v/2}\frako_{K_p}$ then $N\bmx  & b\beta\\ \bar{\beta} & \emx\in p^{v+1}\Z_p$.  So \[
	\calO_p = \qom{b}{\frako_{K_p}}{f\frako_{K_p}}
\]
is a residually inert order with level $p^{v_p(f)+1}$, for $f\in\frako_{K_p}$.

\subsubsection{$K_p$ is unramified}
If $K_p$ is unramified, then a maximal order is given by \[
	\qom{\varpi}{\frako_{K_p}}{\varpi\frako_{K_p}}.
\]
Furthermore, it is well-known that all orders containing the unramified quadratic field extension are isomorphic to $\frako_{K_p}\oplus\varpi^n\frako_{K_p}j$, where $\varpi$ is a uniformizer for $K_p$.  We will represent these orders in this setting via an embedded in the matrix algebra as \[
	\qom{\varpi}{\frako_{K_p}}{\varpi^n\frako_{K_p}}.
\]
Orders of this form have level $p^{2n+1}$.  So consider the order \[
	\calO=\qom{b}{\frako_{K_p}}{f\frako_{K_p}}.
\]
So long as $p\mid b$, this order will have level $p^{2v_p(f)+1}$.

\section{Global order}
Now that we know the form that local orders take, we can examine a global order which has prescribed level locally at each place.  Our global order will be a residually inert order locally for primes $p\mid R_2$, and will be Eichler locally for primes $p\nmid R$, for $a,b\in\Z$ such that $B=\quat{a}{b}{\Q}$ with $a,b\in\Z$ square-free.  Consider the global order $\calO\subset B=\quat{a}{b}{\Q}$ given by \[
	\calO = \calO(a,b,f,g) = \qom{b}{\Z+g\frako_K}{f(\Z+g\frako_K)}.
\]
The localization of this order is \[
	\calO_p = \calO \otimes_{\frako_{\Q}} \frako_{\Q_p} = \calO \otimes \Z_p = \qom{b}{\Z_p+g\frako_{K_p}}{f(\Z_p+g\frako_{K_p})}.
\]
We require $g\in\Z_p^\times$ for the primes where both $K_p$ and $B_p$ are split, so that $\alpha\in \Z_p+\frako_{K_p}=\frako_{K_p}$ and $\beta\in f(\Z_p+\frako_{K_p})=f\frako_{K_p}$.  This yields our order from Section \ref{split-split} with level $p^{v_p(b)+2v_p(f)}$.  Our order $\calO_p$ also has level if $K_p$ is ramified and $B_p$ is split, since its form locally is the same as in Section \ref{split-ram}.

So this order has level $p^{1+2v_p(g)}$ if $p$ is odd or if $p=2$ and $a\equiv 1 \mod 4$, and level $2^{2v_2(g)+2}$ if $p=2$ and $a\equiv 3 \mod 4$, and level $2^{2v_2(g)+3}$ if $p=2$ and $2\mid a$.  Similarly, $\calO_p$ has level if $K_p$ is unramified and $B_p$ is split, giving us level $p^{2v_p(g)}$ if $p$ is odd or $p=2$ with $a\equiv 1 \mod 4$, and $2^{2v_2(g)+2}$ when $p=2$ and $a\equiv 3 \mod 4$, and $2^{2v_2(g)+3}$ when $p=2$ and $2\mid a$.

Locally, our order $\calO_p$ has level if both $K_p$ and $B_p$ are ramified, since our order has the form of the order constructed in \ref{ram-ram}, since $f\in\frako_{K}$ for primes which ramify in both the field and the algebra.  Lastly, requiring $g\in\Z_p^\times$ for the primes where $K_p$ is unramified and $B_p$ is ramified allows $\calO_p$ to have level, since $\alpha\in \Z_p+\frako_{K_p}=\frako_{K_p}$ and $\beta\in f(\Z_p+\frako_{K_p})=f\frako_{K_p}$.  This gives $\calO_p$ level $p^{2v_p(f)+1}$.  We summarize the results from Section \ref{local-orders} via the levels we can achieve locally at each prime, based on the behavior of $K_p$ and $B_p$: \[
	\renewcommand{\arraystretch}{1.3}
		\begin{tabular}{| c || c | c | c | c |}
		\hline
		& $K_p$ split & $K_p$ ramified & $K_p$ unramified\\ \hline
		$B_p$ split & $p^{2v_p(f)+v_p(b)}$ & $p^{2v_p(g)+1}$ & $p^{2v_p(g)}$ \\ \hline
		$B_p$ ramified & $\times$ & $p^{v_p(f)+1}$ & $p^{2v_p(f)+1}$ \\
		\hline
	\end{tabular}
\]
It is important to note the parity that can be achieved in each case; in particular, for $p$ where $K_p$ is unramified and $B_p$ is split, we only obtain even exponents for the local level at $p$.  On the other hand, at $p$ where $K_p$ is ramified and $B_p$ is split, or where $K_p$ is unramified and $B_p$ is ramified, we only obtain odd exponents for the local level at $p$.  When both $K_p$ and $B_p$ are split, we obtain either odd or even exponent, with the parity determined by our selection of $b$ in the representation of $B=\quat{a}{b}{\Q}$.  We have the most freedom when both $K_p$ and $B_p$ are ramified, where we obtain either even or odd exponents, dependent only on the valuation of $f$.

\subsection{Selecting $a,b$}\label{select}

Now consider our quaternion algebra $B$ given via its discriminant $\Delta$, and the level $N$ we desire.  Write $N=R\cdot M$, with primes dividing the discriminant grouped into $R$ and the others into $M$.  Next write $R=R_1\cdot R_2$ and $M=M_1\cdot M_2$, where we group the primes with odd powers into $R_1$ and $M_1$, and the primes with even powers into $R_2$ and $M_2$.  Moreover, we will use $\prod^\prime_{p\mid S}$ to indicate that the product should be taken over all primes $p\in S$ except $p=2$.  We wish to select $a,b\in\Z$ so that (i) $B=\quat{a}{b}{\Q}$ and (ii) primes dividing the level are sorted appropriately into cases which give the correct level that matches the parity of the exponent.

\begin{prop}\label{choose-prop}
Suppose that $B$ is a definite quaternion algebra over $\Q$ with discriminant $\Delta$, and we wish to choose $a,b\in\Z$ so that $B=\quat{a}{b}{\Q}$ and so each prime $p\mid N$ has the appropriate splitting behavior in $K$ and $B$ so that we can achieve local level $p^{v_p(N)}$.  Then we may choose $a,b$ in the following way (noting that $p,q$, and $r$ represent primes) based on our desired behavior at $2$:
	\begin{enumerate}
	
		\item Suppose that $2\nmid\Delta$ and $2$ has an even exponent in $N$ (including the case where $2\nmid N$).  If the product $\prod_{p\mid RM_1} p \equiv 3\mod 4$ then select $a:=-\prod_{p\mid RM_1} p$ and $b:=-q$ with $q$ prime satisfying the conditions 
			\begin{itemize}
				\item $\hs{-q}{p}=-1$ for all $p\mid R$;
				\item $\hs{-q}{p}=1$ for all $p\mid M_1$;
				\item and $q\equiv 1 \mod 8$.
			\end{itemize}
		
			Alternatively, if the product $\prod_{p\mid RM_1} p \equiv 1 \mod 4$ then select $a:=-q\cdot\prod_{p\mid R_2} p$ and $b:=-\prod_{p\mid RM_1} p$ with $q$ prime satisfying the conditions 
			\begin{itemize}
				\item $\hs{q}{p}=(-1)\cdot\hs{-\prod_{r \mid R_2} r}{p}$ for all $p\mid R_1$;
				\item $\hs{q}{p}=(-1)\cdot\hs{-\prod_{r\mid RM_1, r\neq p} r}{p}$ for all $p\mid R_2$;
				\item $\hs{q}{p}=\hs{-\prod_{r\mid R_2} r}{p}$ for all $p\mid M_1$;
				\item If $\prod_{p\mid R_2} p \equiv 1 \mod 4$, then $q\equiv 3 \mod 4$; and if $\prod_{p\mid R_2} p \equiv 3 \mod 4$, then $q\equiv 1 \mod 4$.
			\end{itemize}
		
		\item Suppose that $2\nmid\Delta$ and $2$ has an odd exponent in $N$.  Select $a:=-q\cdot\prod_{p\mid R_2} p$ and $b:=-\prod_{p\mid RM_1} p$ with $q$ prime satisfying the conditions
			\begin{itemize}
				\item $\hs{q}{p}=(-1)\cdot\hs{-\prod_{r\mid R_2} r}{p}$ for all $p\mid R_1$;
				\item $\hs{q}{p}=(-1)\cdot\hs{-\prod_{r\mid RM_1, r\neq p} r}{p}$ for all $p\mid R_2$;
				\item $\hs{q}{p}=\hs{-\prod_{r\mid R_2} r}{p}$ for all $p\mid M_1$, $p\neq 2$;
				\item If $\prod_{p\mid R_2} p \equiv 1 \mod 8$, then $q\equiv 7 \mod 8$; if $\prod_{p\mid R_2} p \equiv 3 \mod 8$, then $q\equiv 5 \mod $; if $\prod_{p\mid R_2} p \equiv 5 \mod 8$, then $q\equiv 3 \mod $; and if $\prod_{p\mid R_2} p \equiv 7 \mod 8$, then $q\equiv 1 \mod 8$.
			\end{itemize}

		\item Suppose that $2\mid\Delta$, $v_2(N)\neq 2$.  Select $a:=-q\cdot\prod_{p\mid R_2} p$ and $b:=-\prod_{p\mid RM_1} p$ with $q$ prime satisfying the conditions 
			\begin{itemize}
				\item $\hs{q}{p}=(-1)\cdot\hs{-\prod^\prime_{r\mid R_2} r}{p}$ for all $p\mid R_1, p\neq 2$;
				\item $\hs{q}{p}=(-1)\cdot\hs{-\prod^\prime_{r\mid RM_1, r\neq p} r}{p}$ for all $p\mid R_2, p\neq 2$;
				\item $\hs{q}{p}=\hs{-\prod^\prime_{r\mid R_2} r}{p}$ for all $p\mid M_1$;
				\item If $v_2(N)=1,3$, choose $q$ so that $a\equiv 5 \mod 8$.  If $v_2(N)>4$ is even, then we have the following for $a'=a/2$ and $b'=b/2$:
				\begin{itemize}
					\item If $b'\equiv 1\mod 8$, then choose $q$ so that $a'\equiv 3$ or $5\mod 8$.
					\item If $b'\equiv 3\mod 8$, then choose $q$ so that $a'\equiv 1$ or $3\mod 8$.
					\item If $b'\equiv 5\mod 8$, then choose $q$ so that $a'\equiv 1$ or $7\mod 8$.
					\item If $b'\equiv 7\mod 8$, then choose $q$ so that $a'\equiv 5$ or $7\mod 8$.
				\end{itemize}
				If $v_2(N)>4$ is odd, then we have the following:
				\begin{itemize}
					\item If $b'\equiv 1\mod 4$, then choose $q$ so that $a\equiv 3\mod 8$.
					\item If $b'\equiv 3\mod 4$, then choose $q$ so that $a\equiv 7\mod 8$.
				\end{itemize}
			\end{itemize}

		\item Lastly, suppose that $2\mid\Delta$ with $v_2(N)=2$.  If $\prod^\prime_{p\mid RM_1} p \equiv 1 \mod 4$, then select $a:=-q\cdot\prod^\prime_{p\mid R_2} p$ and $b:=-\prod^\prime_{p\mid RM_1} p$ with $q$ prime satisfying the conditions 
			\begin{itemize}
				\item $\hs{q}{p}=(-1)\cdot\hs{-\prod^\prime_{r\mid R_2} r}{p}$ for all $p\mid R_1, p\neq 2$;
				\item $\hs{q}{p}=(-1)\cdot\hs{-\prod^\prime_{r\mid RM_1, r\neq p} r}{p}$ for all $p\mid R_2, p\neq 2$;
				\item $\hs{q}{p}=\hs{-\prod^\prime_{r\mid R_2} r}{p}$ for all $p\mid M_1$;
				\item If $\prod^\prime_{p\mid R_2} p \equiv 1 \mod 4$, then $q\equiv 1 \mod 4$; and if $\prod^\prime_{p\mid R_2} p \equiv 3 \mod 4$, then $q\equiv 3 \mod 4$.
			\end{itemize}
			
			Alternatively, if $\prod^\prime_{p\mid RM_1} p \equiv 3 \mod 4$, we cannot construct an order with $v_2(N)=2$.  This is an inherent condition in the structure of the quaternion algebra $B=\quat{a}{b}{\Q}$ and the field $K=\Q(\sqrt{a})$, not specific to our particular construction.
	\end{enumerate}
\end{prop}
	
	Observe that there is a hidden condition that $\hs{b}{q}=1$, so that we obtain disc$\quat{a}{b}{\Q}=\Delta$, which we will verify.  Furthermore, we will observe the behavior of $2$, which will determine its behavior in our algebra.
	
	The conditions on $q$ amount to a finite number of modular congruences, which by Dirichlet's theorem on primes in arithmetic progressions we know have a prime solution $q\neq 2$.  It is worth observing here that if we choose $a$ and $b$ correctly so that the correct prime factors of $a$ and $b$ are ramified (excluding $q$), and if we have the correct ramification or splitting of $2$, we expect that $B_q$ will be split due to the parity of the set of ramified primes.  What we desire in selecting $a$ and $b$ as described above is the following picture of the distribution of primes (ignoring $2$ and $q$): \[
		\renewcommand{\arraystretch}{1.3}
		\begin{tabular}{| c || c | c | c | c |}
		\hline
		& $K_p$ split & $K_p$ ramified & $K_p$ unramified\\ \hline
		$B_p$ split & $p\mid M_2$ & $p\mid M_1$ & $p\mid M_2$ \\
		$B_p$ ramified & $\times$ & $p\mid R$ & $-$ \\
		\hline
		\end{tabular}
	\]
	In the above diagram, the primes dividing $M_2$ are distributed between the $K_p$ split case and the $K_p$ unramified case, since the parity of the exponent we can achieve in those cases is the same.  In particular, according to Proposition \ref{lem-quatsplit} the conditions that $\hs{-q}{p}=-1$ for all $p\mid R$ determine that $B_p$ is ramified for all primes $p\mid R$, and $K_p$ is also ramified for all primes $p\mid R$.  The conditions that $\hs{-q}{p}=1$ for all $p\mid M_1$ determine that $B_p$ is split for all primes $p\mid M_1$, and $K_p$ is ramified for all primes $p\mid M_1$.  These are conditions for all cases except Case 1b, where we will need to compute by hand.  We also need that $B_q$ is split, which is accomplished via quadratic reciprocity: \[
		\hs{a}{q}=\hs{-1}{q}\cdot\prod_{p\mid R} \hs{p}{q} \cdot \prod_{p\mid M_1} \hs{p}{q}
	\]
	\[
		=\hs{-1}{q}\cdot\prod_{p\mid R} (-1)^{\frac{p-1}{2}\frac{q-1}{2}}\hs{q}{p} \cdot \prod_{p\mid M_1} (-1)^{\frac{p-1}{2}\frac{q-1}{2}}\hs{q}{p}
	\]
	\[
		=\hs{-1}{q}\cdot\prod_{p\mid R} (-1)^{\frac{p-1}{2}\frac{q-1}{2}}\hs{-1}{p}\hs{-q}{p} \cdot \prod_{p\mid M_1} (-1)^{\frac{p-1}{2}\frac{q-1}{2}}\hs{-1}{p}\hs{-q}{p}.
	\]
	Notice that we are using quadratic reciprocity assuming that $2\nmid a$.  If $2\mid a$ then we have \[
		\hs{a}{q}=\hs{-2}{q}\cdot\prod_{p\mid R} (-1)^{\frac{p-1}{2}\frac{q-1}{2}}\hs{q}{p} \cdot \prod_{p\mid M_1} (-1)^{\frac{p-1}{2}\frac{q-1}{2}}\hs{q}{p}
	\]
	\[
		=\hs{-2}{q}\cdot\prod_{p\mid R} (-1)^{\frac{p-1}{2}\frac{q-1}{2}}\hs{-1}{p}\hs{-q}{p} \cdot \prod_{p\mid M_1} (-1)^{\frac{p-1}{2}\frac{q-1}{2}}\hs{-1}{p}\hs{-q}{p}.
	\]
	But we have already selected $q$ so that the $\hs{-q}{p}$ in the first product are all $-1$, and the $\hs{-q}{p}$ in the second product are all $1$.  Now there are an odd number of finite ramified primes in $B$, so $\prod_{p\mid R} \hs{q}{p} = -1$.  Thus we have \[
		\hs{a}{q} = \hs{-1}{q}\cdot(-1)\cdot\prod_{p\mid R} (-1)^{\frac{p-1}{2}\frac{q-1}{2}}\hs{-1}{p} \cdot \prod_{p\mid M_1} (-1)^{\frac{p-1}{2}\frac{q-1}{2}}\hs{-1}{p}
	\]
	in the case where $2\nmid a$, and \[
		\hs{a}{q} = \hs{-2}{q}\cdot(-1)\cdot\prod_{p\mid R} (-1)^{\frac{p-1}{2}\frac{q-1}{2}}\hs{-1}{p} \cdot \prod_{p\mid M_1} (-1)^{\frac{p-1}{2}\frac{q-1}{2}}\hs{-1}{p}
	\]
	in the case where $2\mid a$.  In order for $B_q$ to be split we need $\hs{a}{q}=1$, which means that we need \[
		1 = \hs{-1}{q}\cdot(-1)\cdot\prod_{p\mid RM_1} (-1)^{\frac{p-1}{2}\frac{q-1}{2}} \cdot \prod_{p\mid RM_1} \hs{-1}{p}
	\]
	for $2\nmid a$ and \[
		1 = \hs{-2}{q}\cdot(-1)\cdot{\prod_{p\mid RM_1}}^\prime (-1)^{\frac{p-1}{2}\frac{q-1}{2}} \cdot {\prod_{p\mid RM_1}}^\prime \hs{-1}{p}.
	\]
	for $2\mid a$.  These two conditions interact with the behavior of $2$ in $K$ and in $B$, which means that we must consider their behavior together.

\subsubsection{Case 1: $2\nmid\Delta$ and $v_2(N)$ is even}
Suppose we are in Case 1, where we desire $2\nmid\Delta$ and $2$ has an even exponent in $N$.  Furthermore, suppose that the product $\prod_{p\mid RM_1} p \equiv 3 \mod 4$, so that we select $a:=-\prod_{p\mid RM_1} p$ and $b:=-q$, with $q\equiv 1 \mod 8$.  This gives us $a\equiv 1 \mod 4$ and $b\equiv 7 \mod 8$, so $(a,b)_{\Q_2}=(-1)^{\frac{a-1}{2}\frac{b-1}{2}}=1$.  Furthermore, \[
		\hs{a}{q} = \hs{-1}{q}\cdot(-1)\cdot\prod_{p\mid RM_1} (-1)^{\frac{p-1}{2}\frac{q-1}{2}} \cdot \prod_{p\mid RM_1} \hs{-1}{p}
	\]
	\[
		= 1\cdot(-1)\cdot\prod_{p\mid R} \hs{-1}{p} \cdot \prod_{p\mid M_1} \hs{-1}{p} = -1\cdot\prod_{p\mid RM_1} \hs{-1}{p}.
	\]
	Now since $\prod_{p\mid RM_1} p\equiv 3 \mod 4$, an odd number of the $p$ are $\equiv 3 \mod 4$.  Therefore $\prod_{p\mid RM_1} \hs{-1}{p} =-1$, so $\hs{a}{q}=1$ as desired.  Furthermore, $a\equiv 1 \mod 4$ so $2$ is not ramified in $K$.  This gives us \[
		\renewcommand{\arraystretch}{1.3}
		\begin{tabular}{| c || c | c | c | c |}
		\hline
		& $K_p$ split & $K_p$ ramified & $K_p$ unramified\\ \hline
		$B_p$ split & $2, q, p\mid M_2$ & $p\mid M_1$ & $p\mid M_2$ \\
		$B_p$ ramified & $\times$ & $p\mid R$ & $-$ \\
		\hline
	\end{tabular}
	\]
	Alternatively, suppose that the product $\prod_{p\mid RM_1} p \equiv 1 \mod 4$, so that we select $a:=-q\cdot\prod_{p\mid R_2} p$ and $b:=-\prod_{p\mid RM_1} p$.  Then $b\equiv 3 \mod 4$, and we choose $q \mod 4$ so that $a\equiv 1 \mod 4$.  This gives us $(a,b)_{\Q_2}=(-1)^{\frac{a-1}{2}\frac{b-1}{2}}=1$ (so $B_2$ is split) and $K_2$ is not ramified.  Furthermore, we have \[
		\hs{b}{q} = \hs{-\prod_{p\mid RM_1} p}{q} = \hs{-1}{q}\cdot\prod_{p\mid RM_1} \hs{p}{q} = \hs{-1}{q}\cdot\prod_{p\mid R_1} \hs{p}{q} \cdot \prod_{p\mid R_2} \hs{p}{q} \cdot \prod_{p\mid M_1} \hs{p}{q}
	\]
	\[
		= \hs{-1}{q}\cdot \left[\prod_{p\mid RM_1} (-1)^{\frac{p-1}{2}\frac{q-1}{2}}\right] \cdot \left[\prod_{p\mid R_1} \hs{-\prod_{p\mid R_2} p}{p} \cdot (-1)\right]
	\]
	\[
		\cdot \left[\prod_{p\mid R_2} \hs{-\prod_{p\mid R_1M_1} p}{p} \cdot (-1)\right] \cdot \left[\prod_{p\mid M_1} \hs{-\prod_{p\mid R_2} p}{p}\right]
	\]
	via quadratic reciprocity and our assumptions on the values of $\hs{q}{p}$.  Now since $a\equiv 1 \mod 4$, we have $\prod_{p\mid RM_1} (-1)^{\frac{p-1}{2}\frac{q-1}{2}}=1$.  Moreover, we can factor the $-1$'s out of the products above by observing that there are an odd number of primes dividing $R$, so we have an odd number of $-1$'s, giving us a $-1$ factor.  We can also expand the products in the residue symbols above to obtain \[
		\hs{-1}{q}\cdot (-1) \cdot \left[\prod_{r\mid R_1} \hs{-1}{r} \cdot \prod_{p\mid R_2} \hs{p}{r}\right]
	\]
	\[
		\cdot \left[\prod_{r\mid R_2} \hs{-1}{r} \cdot \prod_{p\mid R_1M_1} \hs{p}{r} \right] \cdot \left[\prod_{r\mid M_1} \hs{-1}{r} \cdot \prod_{p\mid R_2} \hs{p}{r}\right].
	\]
	The products above can be written as \[
		\hs{-1}{q}\cdot (-1)\cdot \left[\prod_{p\mid RM_1} \hs{-1}{p}\right]\cdot \left[\prod_{p\mid R_2, r\mid R_1} \hs{p}{r}\right] \cdot \left[\prod_{p\mid R_1, r\mid R_2} \hs{p}{r}\right]
	\]
	\[
		\cdot \left[ \prod_{p\mid M_1, r\mid R_2} \hs{p}{r} \right] \cdot \left[\prod_{p\mid R_2, r\mid M_1} \hs{p}{r}\right].
	\]
	Notice that we have some quadratic reciprocity here; in particular, we have $\hs{p}{r}\cdot\hs{r}{p}$ for all pairs $p\mid R_1$ and $r\mid R_2$, as well as all pairs $p\mid R_2$ and $r\mid M_1$.  Furthermore, since $\prod_{p\mid RM_1} p \equiv 1 \mod 4$, we have $\prod_{p\mid RM_1} \hs{-1}{p} = 1$.  So we obtain \[
		\hs{b}{q} = \hs{-1}{q}\cdot (-1) \cdot \left[\prod_{p\mid R_2, r\mid R_1} (-1)^{\frac{p-1}{2}\frac{r-1}{2}}\right] \cdot \left[\prod_{p\mid R_2, r\mid M_1} (-1)^{\frac{p-1}{2}\frac{r-1}{2}}\right].
	\]
	Now we know that $\prod_{p\mid RM_1} p \equiv 1 \mod 4$, so we have the following cases: 
	\begin{enumerate}
		\item $R_1\equiv 1 \mod 4$, $R_2\equiv 1 \mod 4$, and $M_1\equiv 1 \mod 4$;
		\item $R_1\equiv 1 \mod 4$, $R_2\equiv 3 \mod 4$, and $M_1\equiv 3 \mod 4$;
		\item $R_1\equiv 3 \mod 4$, $R_2\equiv 1 \mod 4$, and $M_1\equiv 3 \mod 4$;
		\item $R_1\equiv 3 \mod 4$, $R_2\equiv 3 \mod 4$, and $M_1\equiv 1 \mod 4$.
	\end{enumerate}
	Furthermore, notice that to obtain $a\equiv 1 \mod 4$, in the first and third cases above we require $q\equiv 3 \mod 4$, and in the second and fourth we require $q\equiv 1 \mod 4$.  Now observe that to obtain $(-1)^{\frac{p-1}{2}\frac{r-1}{2}}=-1$, we need both $p$ and $r$ to be $\equiv 3 \mod 4$.  Moreover, to get $-1$ out of each of the products above, there need to be an odd number of $p\equiv 3 \mod 4$ \textit{and} $r\equiv 3 \mod 4$.  So we can use the above cases to evaluate $\hs{b}{q}$.  In Case 1, $\hs{b}{q} = (-1)(-1)(1)(1) = 1$ as desired; similarly, in Case 2 $\hs{b}{q} = (1)(-1)(1)(-1) = 1$; in Case 3 $\hs{b}{q} = (-1)(-1)(1)(1) = 1$; finally in Case 4 $\hs{b}{q} = (1)(-1)(-1)(1) = 1$.So $B_q$ is split in all cases.  Thus we have the following distribution of primes: \[
		\renewcommand{\arraystretch}{1.3}
		\begin{tabular}{| c || c | c | c | c |}
		\hline
		& $K_p$ split & $K_p$ ramified & $K_p$ unramified\\ \hline
		$B_p$ split & $2, p\mid M_1, M_2$ & $q$ & $p\mid M_2$ \\
		$B_p$ ramified & $\times$ & $p\mid R_2$ & $p\mid R_1$ \\
		\hline
	\end{tabular}
	\]

\subsubsection{Case 2: $2\nmid\Delta$ and $v_2(N)$ is odd}
Suppose now that we are in Case 3, where we desire $2\nmid\Delta$ and $2$ has an odd exponent in $N$.  We select $a:=-q\cdot\prod_{p\mid R_2} p$ and $b:=-\prod_{p\mid RM_1} p$ with $q\mod 8$ so that $a\equiv 1\mod 8$.  Then $2\mid b$.  In this scenario we have \[
		\hs{b}{q} = \hs{-\prod_{p\mid RM_1} p}{q} = \hs{-1}{q}\cdot\prod_{p\mid RM_1} \hs{p}{q} = \hs{-1}{q}\cdot\prod_{p\mid R_1} \hs{p}{q} \cdot \prod_{p\mid R_2} \hs{p}{q} \cdot \prod_{p\mid M_1} \hs{p}{q}
	\]
	This gives us $\hs{b}{q}=1$ as in the previous case.  Moreover, in both cases, since $2\mid b$ we have $2$ non-ramified in $K$.  Now since $2\mid b$, $(a,b)_{\Q_2}=(a,2)_{\Q_2}\cdot(a,b')_{\Q_2}$ for $b'=b/2$.  Now $(a,b')_{\Q_2}=(-1)^{\frac{a-1}{2}\frac{b'-1}{2}}=1$ since $a\equiv 1 \mod 4$, and $(a,2)_{\Q_2}=1$ since $a\equiv 1 \mod 8$.  Thus $B_2$ is split.  So we have the following distribution of primes: \[
		\renewcommand{\arraystretch}{1.3}
		\begin{tabular}{| c || c | c | c | c |}
		\hline
		& $K_p$ split & $K_p$ ramified & $K_p$ unramified\\ \hline
		$B_p$ split & $q, p\mid M_2$ & $2, p\mid M_1$ & $p\mid M_2$ \\
		$B_p$ ramified & $\times$ & $p\mid R_2$ & $R_1$ \\
		\hline
	\end{tabular}
	\]

\subsubsection{Case 3: $2\mid\Delta$ and $v_2(N)$ is odd}
Suppose that we are in Case 3, where we desire $2\mid\Delta$.  We select $a:=-q\cdot\prod_{p\mid R_2} p$ and $b:=-\prod_{p\mid RM_1} p$ with $q\mod 8$ so that $2$ ramifies in $B$ .  Furthermore, we have the requirements $\mod p$ so that the $p\mid R$ are ramified in $B$ and the $p\mid M$ are split in $B$.  Now we have \[
	\hs{b}{q} = \hs{-\prod_{p\mid RM_1} p}{q} = \hs{-1}{q}\cdot\prod_{p\mid RM_1} \hs{p}{q} = \hs{-1}{q}\cdot\prod_{p\mid R_1} \hs{p}{q} \cdot \prod_{p\mid R_2} \hs{p}{q} \cdot \prod_{p\mid M_1} \hs{p}{q}
\]
As before, we obtain $\hs{b}{q} = 1$.  So $B_q$ is split as desired.  Observe that since we have chosen $q\mod 8$ so that $a$ is $3\mod 4$ and $2$ ramifies in $B$, we have $(a,b)_{\Q_2} = -1$ as desired.  In particular, for $v_2(N)=1$ or $3$, we choose $2\mid B$, so our conditions on $q$ make $a\equiv 5\mod 8$ so that $a\equiv 1\mod 4$ and $(a,b)_{\Q_2}=-1$.  For $v_2(N)>4$, we choose $q$ so that $(a,b)_{\Q_2} = (2,2)_{\Q_2}\cdot (a',2)_{\Q_2}\cdot (2,b')_{\Q_2}\cdot (a',b')_{\Q_2} = -1$.  Furthermore, $K_2$ is non-ramified when $2\nmid R_2$, and $K_2$ is ramified when $2\mid R_2$.

\subsubsection{Case 4: $2\mid\Delta$ and $v_2(N)=2$}
It is important to observe that our order at $2$ behaves differently depending on whether $a$ and $b$ are both odd, or $2\mid a$, or $2\mid b$.  The behavior also depends on $a\mod 4$.  In particular, if we desire $v_2(N)=2$, we must select $a$ and $b$ as follows: 

Both $a$ and $b$ must be $\equiv 3\mod 4$, so that $v_2(\disc(K))=2$, and so that $(a,b)_{\Q_2}=(-1)^{\frac{a-1}{2}\frac{b-1}{2}}=-1$.  Furthermore, since we can only choose $a=-\prod^\prime_{p\mid RM_1} p$, $b=-q$ or $a=-q\cdot\prod^\prime_{p\mid R_2} p$, $b=-\prod^\prime_{p\mid RM_1} p$, we therefore need $\prod^\prime_{p\mid RM_1} p \equiv 1\mod 4$.  This is an inherent condition in the structure of the quaternion algebra $B=\quat{a}{b}{\Q}$ and the field $K=\Q(\sqrt{a})$, not specific to our particular construction.  If this product is $\equiv 1\mod 4$, we select $a=-q\cdot\prod^\prime_{p\mid R_2} p$ and $b=-\prod^\prime_{p\mid RM_1} p$ with $q$ so that $a\equiv 3\mod 4$.  Alternatively, if $\prod^\prime_{p\mid RM_1} p\equiv 3\mod 4$, we cannot create an order with $v_2(N)=2$ if $\Delta$ is even.

\subsection{Special Cases}
In some cases, we prefer to restrict our quaternion algebra to simpler scenarios which are common or particularly useful.  In particular, our general construction above can be reduced to two helpful cases: (1) where $\Delta=p$, i.e., where $B$ ramifies at a single prime, and (2) where $K_p$ is unramified wherever $B_p$ is ramified.

\subsubsection{Case 1: $\Delta=p$}\label{sp-case-1}\label{deltap}

Suppose that $\Delta=p$, so our quaternion algebra only ramifies at one place $p$, and furthermore suppose that $N=p^k$.  Then if $p\neq 2$, we can use Case 1 from above to obtain the following:

If $p\equiv 3 \mod 4$ then select $a:=-p$ and $b:=-q$ with $q$ satisfying $\hs{-q}{p}=-1$ and $q\equiv 1 \mod 8$.  If $p\equiv 1 \mod 4$, then choose $a:=-qp$ and $b:=-p$ for even exponents of $p$, and $a:=-q$ and $b:=-p$ for odd exponents of $p$.  In both  scenarios, choose $q$ satisfying $\hs{q}{p} = -1$ and $q\equiv 3 \mod 4$.  In both cases, there is a hidden condition that $\hs{b}{q}=1$ so that $\Delta=p$ as desired.

In the first case, both $2$ and $q$ are split in $K$, and $p$ is ramified in both the field and the algebra.  In the second case, $2$ is split in $K$, $q$ is ramified in $K$ and split in $B$, and $p$ is ramified in $B$ and either ramified or unramified in $K$ (depending on $p \mod 4$).  Then the order \[
		\calO =\left\{\bmx
			\alpha & b\beta\\
			\bar{\beta} & \bar{\alpha}\\
		\emx: \alpha\in\frako_K,\beta\in p^{k-1}\frako_K\right\}
	\]
	has level $qp^k$.  Now the level we have achieved is close to what we desired, but there is an additional $q$.  In a subsequent section we will take care of this and revisit our order.

\subsubsection{Case 2: $K_p$ unramified for $p\mid\Delta$}\label{kp-unram}

Suppose that $R=R_1$, so all primes which ramify in $B$ have odd exponents in the level $N$.  Then we can select $a=-q$ and $b=-\prod_{p\mid RM_1} p$ with the following conditions on $q$:
	\begin{itemize}
		\item $\hs{-q}{p}=-1$ for all $p\mid R$
		\item $\hs{-q}{p}=1$ for all $p\mid M_1$
		\item If $2\mid M_2$, then we require $q\equiv 3 \mod 4$.
		\item If $2\mid M_1$, then we require $q\equiv 7 \mod 8$
		\item If $2\mid R$, then we require $q\equiv 3 \mod 8$.
	\end{itemize}
	
	Conditions 1 and 2 determine that our quaternion algebra is ramified for all $p\mid R$, and split for all $p\mid M_1$, with the exception of $p=2$.  We must verify that $2$ behaves as desired in $B$, and that $B$ splits at $q$.  Now we have \[
		\hs{b}{q} = \hs{-1}{q}\cdot \prod_{p\mid RM_1} \hs{p}{q}
	\]
	if $2\nmid b$, and \[
		\hs{b}{q} = \hs{-2}{q}\cdot {\prod_{p\mid RM_1}}^\prime \hs{-q}{p}\cdot\hs{-1}{p}\cdot(-1)^{\frac{p-1}{2}\frac{q-1}{2}}
	\]
	if $2\mid b$.  In either case, we obtain $\hs{b}{q}=1$ as desired.
	
	In the case where $b$ is odd, $a$ is odd as well, so $(a,b)_{\Q_2}=(-1)^{\frac{a-1}{2}\frac{b-1}{2}}$.  In the first case, $b\equiv 3 \mod 4$, and $a\equiv 1 \mod 4$, so $(-1)^{\frac{a-1}{2}\frac{b-1}{2}}=1$ as desired.  In the second case, $b\equiv 1 \mod 4$, so $(-1)^{\frac{a-1}{2}\frac{b-1}{2}}=1$ as desired.  Thus $2$ is split in $B$ in both cases.  Furthermore, $a\equiv 1 \mod 4$ is required in both cases so that $2$ is not ramified in $K$.  In the other case, where $b$ is even, we use $(a,b)_{\Q_2}=(a,2)_{\Q_2}\cdot (a,b')_{\Q_2}$ (for $b'=b/2$).  Now in both cases, $a\equiv 1\mod 4$, so this yields $(a,b)_{\Q_2}=1$ as desired.  So $2$ is split in $B$.
	
	Now if $2\mid R$, then we have \[
		\hs{b}{q}=\hs{-2}{q}\cdot{\prod_{p\mid RM_1}}^\prime \hs{p}{q} = \hs{-2}{q}\cdot {\prod_{p\mid RM_1}}^\prime \hs{-q}{p}\cdot\hs{-1}{p}\cdot(-1)^{\frac{p-1}{2}\frac{q-1}{2}}
	\]
	\[
		= \hs{-2}{q}\cdot {\prod_{p\mid RM_1}}^\prime \hs{-1}{p}\cdot \left({\prod_{p\mid RM_1}}^\prime (-1)^{\frac{p-1}{2}}\right)^{\frac{q-1}{2}}.
	\]
	Now if $\prod^\prime_{p\mid RM_1} p\equiv 1 \mod 4$, an even number of the $p$ are $\equiv 3 \mod 4$, so $\prod^\prime_{p\mid RM_1} \hs{-1}{p}=1=\prod^\prime_{p\mid RM_1} (-1)^{\frac{p-1}{2}}$, so this reduces to $\hs{b}{q}=\hs{-2}{q}\cdot 1\cdot 1 = \hs{-2}{q}$.  Therefore $q\equiv 3 \mod 8$ is sufficient.
	
	Alternatively, if $\prod^\prime_{p\mid RM_1} p\equiv 3 \mod 4$, an odd number of the $p$ are $\equiv 3 \mod 4$, so $\prod^\prime_{p\mid RM_1} \hs{-1}{p}=-1=\prod^\prime_{p\mid RM_1} (-1)^{\frac{p-1}{2}}$, so this reduces to $\hs{b}{q}=\hs{-2}{q}\cdot -1\cdot (-1)^{\frac{q-1}{2}} = (-1)\cdot \hs{-2}{q}\cdot (-1)^{\frac{q-1}{2}}$.  Therefore $q\equiv 3 \mod 8$ is sufficient.
	
	In both of the cases above, $2\mid b$, so we use $(a,b)_{\Q_2}=(a,2)_{\Q_2}\cdot (a,b')_{\Q_2}$ (for $b'=b/2$).  Now if $2\mid R$, we desire $B_2$ to be ramified.  Now in the first case, $\prod^\prime_{p\mid RM_1} p \equiv 1 \mod 4$, so $b'\equiv 3 \mod 4$.  So $(a,b)_{\Q_2}=(a,2)_{\Q_2}\cdot (a,b')_{\Q_2}=(a,2)_{\Q_2}\cdot (-1)^{\frac{a-1}{2}}$.  Observe that if $q\equiv 3 \mod 8$, then $a\equiv 5 \mod 8$, and either yields $(a,b)_{\Q_2}=-1$ as desired.  So $2$ ramifies in $B$.  In the second case, $\prod^\prime_{p\mid RM_1} p \equiv 3 \mod 4$, so $b'\equiv 1 \mod 4$.  So $(a,b)_{\Q_2}=(a,2)_{\Q_2}\cdot (a,b')_{\Q_2}=(a,2)_{\Q_2}\cdot (-1)^{\frac{a-1}{2}\frac{b-1}{2}}=(a,2)_{\Q_2}$.  Observe that if $q\equiv 3 \mod 8$, then $a\equiv 5 \mod 8$, and either yields $(a,b)_{\Q_2}=-1$ as desired.  So $2$ ramifies in $B$.  Furthermore, we desire that $2$ is unramified in $K$, so we need $q\equiv 3 \mod 4$.  Therefore, we need $q\equiv 3 \mod 8$.
	
	Now we have selected $a,b$ so that $B=\quat{a}{b}{\Q}$, So we have the following distribution of primes: \[
		\renewcommand{\arraystretch}{1.3}
		\begin{tabular}{| c || c | c | c | c |}
		\hline
		& $K_p$ split & $K_p$ ramified & $K_p$ unramified\\ \hline
		$B_p$ split & $p\mid M_1, p\mid M_2$ & $q$ & $p\mid M_2$ \\
		$B_p$ ramified & $\times$ & $-$ & $p\mid R$ \\
		\hline
	\end{tabular}
	\]

	Now the order \[
		\calO =\left\{\bmx
			\alpha & b\beta\\
			\bar{\beta} & \bar{\alpha}\\
		\emx: \alpha\in\Z+g\frako_K,\beta\in f(\Z+g\frako_K)\right\}
	\]
	has level $qN$ for \[
		f=\epsilon\cdot{\prod_{p\mid R}}^\prime p^{\frac{v_p(N)-1}{2}}\cdot\prod_{\substack{p\mid M_2,\\ \hs{a}{p}=1}} p^{v_p(N)/2}, \quad 
		g=\prod_{\substack{p\mid M_2,\\ \hs{a}{p}=-1}} p^{v_p(N)/2}\cdot \prod_{p\mid M_1} p^{\frac{v_p(N)-1}{2}},
	\]
	with \[
		\epsilon = \begin{cases}
			2^{\frac{v_2(N)-3}{2}} & \text{if } v_2(N)\geq 3 \text{ odd}\\
			2^{\frac{v_2(N)}{2}-2} & \text{if } v_2(N)\geq 3 \text{ even}\\
			1 & \text{else}.
		\end{cases}
	\]
	Now the level we have achieved is close to what we desired, but there is a $q$ factor separating us from our ultimate goal.  In a subsequent section we will take care of this and revisit our order.

\subsection{Basis for $\calO$}
We know that $\frako_K$ has basis $1,\sqrt{a}$ if $a\equiv 3 \mod 4$, and $1,\frac{1+\sqrt{a}}{2}$ if $a\equiv 1 \mod 4$.  In the first case, $\Z+g\frako_K$ has basis $1,g\sqrt{a}$ and $f(\Z+g\frako_K)$ has basis $f, fg\sqrt{a}$; in the second case, $\Z+g\frako_K$ has basis $1,g\left(\frac{1+\sqrt{a}}{2}\right)$ and $f(\Z+g\frako_K)$ has basis $f,fg\left(\frac{1+\sqrt{a}}{2}\right)$.  So a basis for \[
	\calO =\left\{\bmx
			\alpha & b\beta\\
			\bar{\beta} & \bar{\alpha}\\
		\emx: \alpha\in\Z+g\frako_K,\beta\in f(\Z+g\frako_K)\right\}
\]
is \[
	\bmx 1 & \\ & 1\emx, \bmx g\sqrt{a} & \\ & -g\sqrt{a}\emx, \bmx & fb\\ f & \emx, \text{ and } \bmx & fgb\sqrt{a}\\ -fg\sqrt{a} & \emx
\]
in the case where $a\equiv 3 \mod 4$, and \[
	\bmx 1 & \\ & 1\emx, \bmx g\left(\frac{1+\sqrt{a}}{2}\right) & \\ & g\left(\frac{1-\sqrt{a}}{2}\right) \emx, \bmx & fb\\ f & \emx, \text{ and } \bmx & fgb\left(\frac{1+\sqrt{a}}{2}\right)\\ fg\left(\frac{1-\sqrt{a}}{2}\right) & \emx
\]
in the case where $a\equiv 1 \mod 4$.  Using the fact that \[
	\bmx \sqrt{a} & \\ & -\sqrt{a} \emx \cdot \bmx & b\\ 1 & \emx = \bmx & b\sqrt{a}\\ -\sqrt{a} & \emx
\]
we can convert these elements from their matrix representations to the standard representations using $1, i, j,$ and $k$.  We also observe that since $f\in\frako_K$, then we will obtain powers of the factors of $a$, encoded by $h$ (a square-free integer).  Notice that locally when $p$ is ramified in $K$--so that $\calO_p$ is a residually inert order--we have basis \[
	\calO_p = \Z\la 1,i,p^{v_p(h)}\cdot fj, fk \ra,
\]
with level $p^{2v_p(f)+2}$ when $p$ is ramified in $B$.  This allows us to simplify our description of $f$ to be in $\Z$ rather than in $\frako_K$, by incorporating a third, square-free integer $h$ to do globally what the $p$ did locally here.  So, if $a\equiv 3 \mod 4$, our basis becomes \[
	\calO=\Z\la1, gi, fhj, fgk\ra.
\]
Alternatively, if $a\equiv 1 \mod 4$, our basis becomes \[
	\calO=\Z\left\la1, g\left(\frac{1+i}{2}\right), fhj, fg\left(\frac{hj+k}{2}\right)\right\ra.
\]
We can calculate the discriminant of $\calO_p$ by observing that if $p\neq 2$ we have $\calO_p=\Z_q\la 1, gi, fgh\cdot j, fk\ra$ and $\text{disc}(\calO_p)=\sqrt{\det{(\alpha_i\alpha_j)}}=ab\cdot f^2g^2h$.  If $p=2$, then $\disc(\calO_2)=4$ if $a\equiv 3 \mod 4$, and $\disc(\calO_2)=1$ if $a\equiv 1 \mod 4$.  If $2\mid a$, then $\disc(\calO_2)=8$.

\section{Lowering the Level}
In Section \ref{select} you may notice that the prime $q$ was used to manipulate our quaternion algebra $\quat{a}{b}{\Q}$ so that we obtained the discriminant $\Delta$ as we desired, as well as distributing the primes in $R_1$, $R_2$, $M_1$ and $M_2$ properly to obtain the desired parity for each exponent.  Furthermore, $q$ was often used to manipulate the behavior of $2$ in both the quadratic field $K(\sqrt{a})$ and $B$ to obtain the desired behavior of $2$ in the level.  The selection of $q$ played a central role in achieving these results, yet as you may note from the distributions of primes in each of the cases from the previous section the behavior of $K_q$ and $B_q$ is such that while $B_q$ is always split we can obtain at minimum a level of $q$ for the localization of our order $\calO_q$.  So to obtain our desired level $N$, we must lower the level of $\calO_q$ from $q$ to $1$, which in turn lowers the level of $\calO$ from $qN$ to $N$.  Observe that since we have chosen $q\neq 2$, $\calO_q$ has basis $1,i,j,k$ in all cases.  This allows us to apply a technique from Voight (\cite{voight-max}) to find a maximal order $\calO_q'$ containing $\calO_q$.  From \cite{voight-max}, Algorithm 7.10, we will compute a $q$-maximal order containing $\calO$ by adjoining a special element to our order.  Now in Cases 1a, 2, and 3, $q\mid b$, and in Cases 1b $q\mid a$.  Now since $q$ is odd, we are in Step 2 of the algorithm, where we swap $i$ for $j$ or $k$ so that $\text{ord}_q(a)=0$.  So in Cases 1a, 2, and 3 $q\nmid a$, so we do not need to swap anything.  If we are in Cases 1b, we swap $i$ for $j$ locally, which globally swaps $gi$ and $fj$.  In both of these cases, $\text{ord}_q(b)=1$ and (after swapping if necessary) $\hs{a}{q}=1$, so next we solve $x^2\equiv a \mod q$ for $x\in \Z/q\Z$ and adjoin $q^{-1}(x-i)j$ locally.  In order to adjoin this element globally to $\calO$ without altering $\calO_p$ ($p\neq q$) we adjoin $fgq^{-1}(x-i)j$ globally in Cases 1a, 2, and 3, and we adjoin $fgq^{-1}(x-j)i$ globally in Cases 1b.

Adjoining $fgq^{-1}(x-i)j$ globally to our order does not affect the order at places $p\neq q$, since $q\in\Z_p^\times$ and $x\in\Z_p$, so we have \[
	fgq^{-1}(x-i)j=fgq^{-1}(xj-k)=q^{-1}(xg(fj)-fgk)\in\calO_p
\]
for all $p\neq q$.  Similarly, adjoining $fgq^{-1}(x-j)i$ globally to our order does not affect the order at places $p\neq q$ since \[
	fgq^{-1}(x-j)i=fgq^{-1}(xi+k)=q^{-1}(xf(gi)+fgk)\in\calO_p
\]
for all $p\neq q$.  In particular, \[
	\Z_p\la 1, gi, fhj, fgq^{-1}(x-i)j\ra = \Z_p\la 1, gi, fhj, fg(xj-k)\ra = \Z_p\la 1, gi, fhj, gx(fj)-fgk\ra
\]
\[
	= \Z_p\la 1, gi, fhj, fgk\ra = \calO_p
\]
and \[
	\Z_p\la 1, gi, fhj, fgq^{-1}(x-j)i\ra = \Z_p\la 1, gi, fhj, fg(xi+k)\ra = \Z_p\la 1, gi, fhj, fx(gi)+fgk\ra
\]
\[
	= \Z_p\la 1, gi, fhj, fgk\ra = \calO_p.
\]
Therefore we have $\calO=\Z\la 1, gi, fhj, fgk \ra \subset \calO'=\Z\la 1,gi,fhj,fgq^{-1}(x-i)j\ra$.  On the other hand, away from $q$, $\calO'_p=\calO_p$; at $q$, $\calO_q\subset\calO'_q$ with $\calO_q'$ maximal as desired, since \[
	\Z_q\la 1, gi, fhj, fgq^{-1}(x-i)j\ra = \Z_q\la 1, i, j, q^{-1}(x-i)j\ra
\]
and \[
	\Z_q\la 1, gi, fhj, fgq^{-1}(x-j)i\ra = \Z_q\la 1, i, j, q^{-1}(x-j)i\ra.
\]
So since $\calO'$ is unchanged from $\calO$ for $p\neq q$, while $\calO'$ has level 1 at $q$ whereas $\calO$ has level $q$, we have $\lev(\calO)=R_1R_2\cdot M_1M_2=N$ as desired.

So globally if $a\equiv 3 \mod 4$ we can compute our basis as \[
	\calO' = \begin{cases}
		\Z\la 1, gi, \frac{fgh(j+xk)}{q}, fk \ra & \text{ if } q\mid b\\
		\Z\la 1, \frac{gh(i+tk)}{q}, fgj, fk \ra & \text{ if } q\mid a\\
	\end{cases},
\]
where $t\in\Z$ comes from the Euclidean algorithm for writing 
	\begin{equation}\label{t}
	1 = s(q)+t(hx).
	\end{equation}
However, if $a\equiv 1 \mod 4$, so $\calO$ has basis \[
	\calO = \Z\left\la 1, g\cdot\frac{1+i}{2}, fhj, fg\cdot\frac{hj+k}{2} \right\ra,
\]
we must calculate the basis obtained by adjoining $q^{-1}(x-i)j$ or $q^{-1}(x-j)i$ to $\calO$ by using the Hermite normal form.  The Hermite normal forms give us a basis for our order: \[
	\calO' = \begin{cases}
		\Z\left\la 1, \frac{g(1+i)}{2}, \frac{fgh(j+uk)}{2q}, fk \right\ra & \text{if } q\mid b\\
		\Z\left\la \frac{q+gi+2gzk}{2q}, \frac{g(2i+z'k)}{2q}, \frac{f(hj+gk)}{2}, fgk \right\ra & \text{if } q\mid a\\
	\end{cases},
\]
where $u$ is given by using the Euclidean algorithm to write $v(q)+w(2x)=1$, and setting $0\leq u < 2q$ such that 
	\begin{equation}\label{u}
		u\equiv vq+2w \mod 2q,
	\end{equation}
$z$ is given by the Euclidean algorithm for writing 
	\begin{equation}\label{y}
		y(-q)+z(2x)=1,
	\end{equation}
and where $z'$ is given by choosing $0\leq z' < 2q$ with
	\begin{equation}\label{y'}
		z'\equiv 4z \mod 2q.
	\end{equation}

\section{Main result}
In Sections 4 we constructed orders of level $qN$, and in the previous section we lowered the level of our order at $q$ from $q$ to $1$.  Therefore we have the following: 

\begin{thm}\label{big-thm}
Select $a,b$ to represent our quaternion algebra as stated in Proposition \ref{choose-prop}, and put \[
		f=\epsilon\cdot{\prod_{p\mid R_1}}^\prime p^{\frac{v_p(N)-1}{2}}\cdot{\prod_{p\mid R_2}}^\prime p^{v_p(N)/2-1}\cdot\prod_{\substack{p\mid M_2,\\ \hs{a}{p}=1}} p^{v_p(N)/2}, \quad 
		g=\prod_{\substack{p\mid M_2,\\ \hs{a}{p}=-1}} p^{v_p(N)/2}\cdot \prod_{p\mid M_1} p^{\frac{v_p(N)-1}{2}}
	\]
and \[
	h = \begin{cases} \ds\prod_{p\mid R_2} p^{1-v_p(b)} & \text{if } v_2(N)\neq 2\\ {\ds\prod_{p\mid R_2}}^\prime p^{1-v_p(b)} & \text{if } v_2(N)= 2\\ \end{cases} \text{ with } \epsilon = \begin{cases}
		2^{\frac{v_2(N)-3}{2}} & \text{if } v_2(N)\geq 3 \text{ odd}\\
		2^{\frac{v_2(N)}{2}-2} & \text{if } v_2(N)\geq 3 \text{ even}\\
		1 & \text{else}
	\end{cases}.
\]

and select $x\in\Z$ with $x^2\equiv a \mod q$ if $q\mid b$, and $x^2\equiv b \mod q$ if $q\mid a$.  Then the order \[
	\calO = \begin{cases}
		\Z\left\la \frac{q+gi+2gzk}{2q}, \frac{g(2i+z'k)}{2q}, \frac{f(hj+gk)}{2}, fgk \right\ra & \text{ if } q\mid a \text{ and } a\equiv 1\mod 4\\
		\Z\la 1, \frac{gh(i+tk)}{q}, fgj, fk \ra & \text{ if } q\mid a \text{ and } a\equiv 3\mod 4\\
		\Z\left\la 1, \frac{g(1+i)}{2}, \frac{fgh(j+uk)}{2q}, fk \right\ra & \text{ if } q\mid b \text{ and } a\equiv 1\mod 4\\
		\Z\la 1, gi, \frac{fgh(j+xk)}{q}, fk \ra & \text{ if } q\mid b \text{ and } a\equiv 3\mod 4\\
	\end{cases}
\]
has level $N$ in $B$, with $u$ given by \eqref{u}, $z$ given by \eqref{y}, $z'$ given by \eqref{y'}, and $t$ given by \eqref{t}.
\end{thm}

Theorem \ref{thm-1-1} follows from descending Theorem \ref{big-thm} to $\Delta=p$.

When $\Delta\mid R_1$--i.e. when $K_p$ is unramified for all $p\mid\Delta$--Theorem \ref{big-thm} descends to 
\begin{thm}
When we select $a,b$ to represent our quaternion algebra as stated in Section \ref{kp-unram}, we choose \[
		f=\epsilon\cdot\prod_{p\mid R} p^{\frac{v_p(N)-1}{2}}\cdot\prod_{\substack{p\mid M_2,\\ \hs{a}{p}=1}} p^{v_p(N)/2}, \quad 
		g=\prod_{\substack{p\mid M_2,\\ \hs{a}{p}=-1}} p^{v_p(N)/2}\cdot \prod_{p\mid M_1} p^{\frac{v_p(N)-1}{2}}
	\]
with \[
	\epsilon = \begin{cases}
		2^{\frac{v_2(N)-3}{2}} & \text{if } v_2(N)\geq 3 \text{ odd}\\
		2^{\frac{v_2(N)}{2}-2} & \text{if } v_2(N)\geq 3 \text{ even}\\
		1 & \text{else}
	\end{cases},
\]
and select $x\in\Z$ with $x^2\equiv a \mod q$ if $q\mid b$, and $x^2\equiv b \mod q$ is $q\mid a$.  Then the order \[
	\calO = \begin{cases}
		\Z\left\la \frac{q+gi+2gzk}{2q}, \frac{g(2i+z'k)}{2q}, \frac{f(j+gk)}{2}, fgk \right\ra & \text{ if } a\equiv 1\mod 4\\
		\Z\la 1, \frac{g(i+tk)}{q}, fgj, fk \ra & \text{ if } a\equiv 3\mod 4\\
	\end{cases}
\]
has level $N$ in $B$, with $z$ given by \eqref{y}, $z'$ given by \eqref{y'}, and $t$ given by \eqref{t}.
\end{thm}

These results have been checked for $\Delta<1000$ and $N<10,000$ by constructing the order prescribed above in Sage and computing its discriminant, matching it to the level $N$.  Note that I have provided the code for the general construction of an order with level $L$ detailed in my result, available at \url{http://math.ou.edu/~jwiebe/}.

\section{Examples}
Now that we have our order $\calO$ of level $N$, we can use it to construct spaces of modular forms of level $N$ using Brandt matrices (or theta series); see \cite{pizer-algorithm}, \cite{hps-orders} when $B$ has prime discriminant, and see \cite{martin-basis} for arbitrary $B$.  Note that there are other approaches to this, including a technique of Dembele \cite{dembele} which only requires the use of maximal orders.  However, our result also allows us to compute quaternionic modular forms via Brandt matrices, and also solves the quaternionic analog of the classical problem of finding bases for orders in number fields.

We conclude by presenting a couple of examples of finding bases
of orders, and indicate how this is used to compute modular forms of matching level.

\begin{ex}[$\Delta=3$ and $N=27$]\end{ex}
Suppose that $\Delta=3$ and $N=27$.  We can compute the class number (see \cite{hps-memoirs}), and obtain $H=2$.  We can compute $a,b$ and $\calO$ using the case outlined in Section \ref{sp-case-1} to obtain $a=-7$, $b=-3$ and use Theorem \ref{big-thm} to obtain \[
	\calO = \Z\left\la \frac{1+i}{2}, i, \frac{3j+309k}{146}, 3k \right\ra.
\]
Using Magma we obtain the following via \texttt{M:=BrandtModule(O)} and \texttt{HeckeOperator(M,p)}: \[
	T_1 = \bmx 1 & 0 & 0 & 0\\ 0 & 1 & 0 & 0\\ 0 & 0 & 1 & 0\\ 0 & 0 & 0 & 1\\ \emx, \hspace{10pt} 
	T_2 = \bmx 1 & 2 & 1 & 2\\ 1 & 2 & 1 & 2\\ 1 & 2 & 1 & 2\\ 1 & 2 & 1 & 2\\ \emx, \hspace{10pt} 
	T_3 = \bmx 0 & 0 & 0 & 0\\ 0 & 0 & 0 & 0\\ 0 & 0 & 0 & 0\\ 0 & 0 & 0 & 0\\ \emx, \hspace{10pt} 
	T_4 = \bmx 4 & 12 & 6 & 12\\ 6 & 10 & 6 & 12\\ 6 & 12 & 4 & 12\\ 6 & 12 & 6 & 10\\ \emx,
\]
\[
	T_5 = \bmx 2 & 4 & 2 & 4\\ 2 & 4 & 2 & 4\\ 2 & 4 & 2 & 4\\ 2 & 4 & 2 & 4\\ \emx, \hspace{10pt} 
	T_6 = \bmx 0 & 0 & 0 & 0\\ 0 & 0 & 0 & 0\\ 0 & 0 & 0 & 0\\ 0 & 0 & 0 & 0\\ \emx, \hspace{10pt} 
	T_7 = \bmx 2 & 6 & 2 & 6\\ 3 & 5 & 3 & 5\\ 2 & 6 & 2 & 6\\ 3 & 5 & 3 & 5\\ \emx, \hspace{10pt} 
	T_8 = \bmx 32 & 64 & 32 & 64\\ 32 & 64 & 32 & 64\\ 32 & 64 & 32 & 64\\ 32 & 64 & 32 & 64\\ \emx, \dots
\]
which yield the Eisenstein series with $a_p=p+1$ for $p\neq 3$, as well as the cusp form \[
	f = q - 2q^4 - q^7 + 5q^{13} + \dots.
\]
These are both modular forms of weight $2$ and level $27$, whose $p$th Fourier coefficient is an eigenvalue of the Hecke operator $T_p$ above ($p\neq 3$).

\begin{ex}[$\Delta=7$ and $N=49$]\end{ex}
Suppose that $\Delta=7$ and $N=49$.  We can compute $a,b$ and $\calO$ using (again) the case outlined in Section \ref{sp-case-1} to obtain $a=-7$, $b=-11$ and use Theorem \ref{big-thm} to obtain \[
	\calO = \Z\left\la 1,\frac{1+i}{2},\frac{7(j-5k)}{22},k \right\ra.
\]
Using Magma as in the previous example gives us \[
	T_1 = \bmx 1 & 0 & 0 & 0\\ 0 & 1 & 0 & 0\\ 0 & 0 & 1 & 0\\ 0 & 0 & 0 & 1\\ \emx, \hspace{10pt}
	T_2 = \bmx 2 & 1 & 0 & 0\\ 1 & 2 & 0 & 0\\ 0 & 0 & 2 & 1\\ 0 & 0 & 1 & 2\\ \emx, \hspace{10pt} 
	T_3 = \bmx 0 & 0 & 2 & 2\\ 0 & 0 & 2 & 2\\ 2 & 2 & 0 & 0\\ 2 & 2 & 0 & 0\\ \emx, \hspace{10pt} 
	T_4 = \bmx 3 & 4 & 0 & 0\\ 4 & 3 & 0 & 0\\ 0 & 0 & 3 & 4\\ 0 & 0 & 4 & 3\\ \emx,
\]
\[
	T_5 = \bmx 0 & 0 & 3 & 3\\ 0 & 0 & 3 & 3\\ 3 & 3 & 0 & 0\\ 3 & 3 & 0 & 0\\ \emx, \hspace{10pt} 
	T_6 = \bmx 0 & 0 & 6 & 6\\ 0 & 0 & 6 & 6\\ 6 & 6 & 0 & 0\\ 6 & 6 & 0 & 0\\ \emx, \hspace{10pt} 
	T_7 = \bmx 2 & 2 & 2 & 2\\ 2 & 2 & 2 & 2\\ 2 & 2 & 2 & 2\\ 2 & 2 & 2 & 2\\ \emx, \hspace{10pt}
	T_8 = \bmx 6 & 9 & 0 & 0\\ 9 & 6 & 0 & 0\\ 0 & 0 & 6 & 9\\ 0 & 0 & 9 & 6\\ \emx, \dots
\]
which yield the Eisenstein series with $a_p=p+1$ for $p\neq 7$, as well as the cusp form \[
	f = q  + q^2 - q^4 - 3q^8 - 3q^9 + 4q^{11} + \dots
\]
These are both modular forms of weight $2$ and level $49$, whose $p$th Fourier coefficient is an eigenvalue of the Hecke operator $T_p$ above ($p\neq 7$).

\begin{ex}[$\Delta=70$ and $N=2\cdot5^2\cdot7^5\cdot11\cdot 23^2$]\end{ex}
Suppose that $\Delta=70$ and $N=2\cdot5^2\cdot7^5\cdot11\cdot 23^2$.  Since $2\mid\Delta$ and $v_2(N)=1$, we are in Case 3, where we select $a=-q\cdot\prod_{p\mid R_2} p$ and $b=-\prod_{p\mid RM_1} p$.  Our conditions on $q$ we compute as:

\begin{enumerate}
	\item $\hs{q}{7} = (-1)\cdot \hs{-5}{7}$; and 
	\item $\hs{q}{5} = (-1)\cdot \hs{-7\cdot 11}{5}$; and 
	\item $\hs{q}{11} = \hs{-5}{11}$.
\end{enumerate}

We also choose $q\equiv 7\mod 8$ so that $a\equiv 5\mod 8$.  So this gives us a set of congruences where $q$ is nonsquare $\mod 7$, a square $\mod 5$, and nonsquare $\mod 11$.  So if $q\equiv 7\mod 8$, $q\equiv 3 \mod 7$, $q\equiv 2\mod 5$, and $q\equiv 2\mod 11$, we obtain $q=1487$.  So $a=-1487\cdot 5$ and $b=-2\cdot 5\cdot 7\cdot 11$.  This gives us $B = \quat{-7435}{-770}{\Q}$ with $\Delta=70$ as desired.  Then $f=23\cdot 7^2$, $g=1$, and $h=1$.  Now we need to find $x$ so that $x^2\equiv -770\mod 1487$, which gives us $x=593$.  Next we use the extended Euclidean algorithm to compute $d=y(-1487)+z(2\cdot 593)$, which gives us $z=1156$.  Then $c=2z=2312$.  So now we can construct our order: \[
	\calO = \Z\left\la \frac{1487+i+2\cdot578k}{2974}, \frac{i+ 1156k}{1487}, \frac{1127j + 1127k}{2}, 1127k \right\ra.
\]
This order has level $N=4889996650=2\cdot5^2\cdot7^5\cdot11\cdot23^2$ as desired.  We can continue in the process illustrated in the previous examples, using this order, employing Magma to compute the Hecke operators for this order, and compute from them the modular forms of weight $2$ and level $N$.\\

\pagebreak

\bibliographystyle{unsrt}
\bibliography{nonmaximal-orders}

\end{document}